\documentclass{gtart}    

\def\ifplaintex{\expandafter\ifx\csname documentclass\endcsname\relax}

\ifplaintex 
\else
\fi

\def\gtm{{\mathsurround=0pt\it $\cal G\mskip-2mu$eometry \&\ 
$\cal T\!\!$opology $\cal M\mskip-1mu$onographs}}    

\def\gtp{{\mathsurround=0pt\it $\cal G\mskip-2mu$eometry \&\ 
$\cal T\!\!$opology $\cal P\!$ublications}}  

\def\recd{{\small Received:\qua\receiveddate\ifx\reviseddate\relax
\else\qquad Revised:\qua\reviseddate\fi\par}} 

\def\volumenumber#1{\def\thevolumenumber{#1}}
\def\volumeyear#1{\def\thevolumeyear{#1}}
\def\volumename#1{\def\thevolumename{#1}}
\def\papernumber#1{\def\thepapernumber{#1}}
\def\pagenumbers#1#2{\def\startpage{#1}\def\finishpage{#2}}
\def\published#1{\def\publishdate{#1}}
\def\received#1{\def\receiveddate{#1}}

\def\accepted#1{\def\accepteddate{#1}}

\def\asciiaddress#1{\def\theasciiaddress{#1}}

\let\\\par
\let\thevolumenumber\relax\let\thepapernumber\relax
\let\thevolumeyear\relax\let\startpage\relax
\let\finishpage\relax\let\publishdate\relax\let\receiveddate\relax
\let\reviseddate\relax\let\accepteddate\relax\let\theasciititle\relax
\let\theasciiauthors\relax\let\theasciiaddress\relax
\let\theasciiabstract\relax

\let\theerratum\relax\let\theasciiemail\relax
\let\theshortauthors\relax\let\theshorttitle\relax

\def\startpage{1}\def\finishpage{15}\def\thepapernumber{77}

\volumenumber{2}
\volumename{Proceedings of the Kirbyfest}
\volumeyear{1999}

\long\def\maketitlep{   

\count0=\startpage

\gtm\nl        
{\small Volume \thevolumenumber: \thevolumename\nl 
\ifx\theerratum\relax\else Erratum \erratumnumber\nl\fi
Pages \startpage--\finishpage\nl}

\vglue 0.1truein   

{\parskip=0pt\leftskip 0pt plus 1fil\def\\{\par\smallskip}{\ifplaintex\large
\else\Large\fi\bf\thetitle}\par\medskip}   
\vglue 0.05truein 

{\parskip=0pt\leftskip 0pt plus 1fil\def\\{\par}{\sc\theauthors}
\par\medskip}%
 
\vglue 0.03truein 

{\small\leftskip 25pt\rightskip 25pt{\bf Abstract}\stdspace\theabstract

{\bf AMS Classification}\stdspace\theprimaryclass
\ifx\thesecondaryclass\relax\else; \thesecondaryclass\fi\par
{\bf Keywords}\stdspace \thekeywords\par}\vglue 7pt

}   

\font\phead=cmsl9 scaled 950
\font=cmsl9 scaled 1050
\font\pnum=cmbx10 scaled 913
\font\lnum=cmbx10 
\font\pfoot=cmsl9 scaled 950
\font=cmsl9 scaled 1050
\ifplaintex
\headline{\vbox to 0pt{\vskip -4.5mm\line{\small\phead\ifnum
\count0=\startpage ISSN 1464-8997 (on line)
1464-8989 (printed) \hfill {\pnum\folio}\else\ifodd\count0\def\\{ }%
\ifx\theshorttitle\relax\thetitle\else\theshorttitle\fi\hfill{\pnum\folio}
\else\def\\{ and }{\pnum\folio}\hfill\ifx\theshortauthors\relax\theauthors
\else\theshortauthors\fi\fi\fi}\vss}}
\footline{\vbox to 0pt{\vglue 0mm\line{\small\pfoot\ifnum\count0=\startpage
Published \publishdate:\qua\copyright\ \gtp\hfill\else
\gtm, Volume \thevolumenumber\ (\thevolumeyear)\hfill\fi}\vss
}}
\else
\makeatletter
\def\@oddhead{{\small\lhead\ifnum\count0=\startpage ISSN 1464-8997 (on line)
1464-8989 (printed) \hfill {\lnum\number\count0}\else\ifodd\count0
\def\\{ }\ifx\theshorttitle\relax \thetitle \else\theshorttitle\fi\hfill
{\lnum\number\count0}\else\def\\{ and }{\lnum\number\count0}
\hfill\ifx\theshortauthors\relax 
\theauthors\else\theshortauthors\fi\fi\fi}}\def\@evenhead{@oddhead}
\def\@oddfoot{\small\lfoot\ifnum\count0=\startpage Published \publishdate:\qua\copyright\ \gtp\hfill\else
\gtm, Volume \thevolumenumber\ (\thevolumeyear)\hfill\fi}
\def\@evenfoot{@oddfoot}
\makeatother
\fi

\let\maketitlepage\maketitlep
\let\makeshorttitle\maketitlepage
\let\maketitle\maketitlepage

\newwrite\gtoutfile
\long\gdef\makeheadfile{  
{\def\\{, }\def\s{ }
\immediate\openout\gtoutfile head.xxx
\immediate\write\gtoutfile{To: math@arxiv.org}
\immediate\write\gtoutfile{Subject: put OR rep NNNNN:ppppp}
\immediate\write\gtoutfile{--text follows this line--}
\immediate\write\gtoutfile{Proxy-for: \ifx\theasciiauthors\relax
\theauthors\else\theasciiauthors\fi\s<\ifx\theasciiemail\relax\theemail\else\theasciiemail\fi>}
\immediate\write\gtoutfile{\noexpand\\}
\immediate\write\gtoutfile{Authors: \ifx\theasciiauthors\relax
\theauthors\else\theasciiauthors\fi}
{\def\\{ }\immediate\write\gtoutfile{Title: \ifx\theasciititle\relax
\thetitle\else\theasciititle\fi}}
\immediate\write\gtoutfile{Subj-class: GT or SG, GR etc}
\immediate\write\gtoutfile{MSC-class: \theprimaryclass\ifx\thesecondaryclass\relax\else, \thesecondaryclass\fi}
\immediate\write\gtoutfile{Journal-ref: Geom. Topol. Monogr. \thevolumenumber\s
(\thevolumeyear) \startpage-\finishpage}
\immediate\write\gtoutfile{Comments: Published by Geometry and Topology Monographs at}
\immediate\write\gtoutfile{\s\s\s  http://www.maths.warwick.ac.uk/gt/GTMon\thevolumenumber/paper\thepapernumber.abs.html}
\immediate\write\gtoutfile{\noexpand\\}
\immediate\write\gtoutfile{}
\ifx\theasciiabstract\relax
\immediate\write\gtoutfile{\theabstract}\else
\immediate\write\gtoutfile{\theasciiabstract}\fi
\immediate\write\gtoutfile{}
\immediate\write\gtoutfile{\noexpand\\}
\immediate\write\gtoutfile{}
\immediate\closeout\gtoutfile}}  

\def\maketitlepage{\maketitlep\makeheadfile}
\let\makeshorttitle\maketitlepage
\let\maketitle\maketitlepage

\volumenumber{4}
\volumename{Invariants of knots and 3-manifolds (Kyoto 2001)}
\volumeyear{2002}
\papernumber{19} 
\pagenumbers{295}{302} 
\received{27 November 2001}
\accepted{22 July 2002}
\published{13 October 2002}

\def\card{{  \rm {card}}}

\def\sh#1{\section*{#1}}

\def\Eul{{ \rm  {Eul}}}

\def\Pf {{ \rm  {Pf}}}

\def\sign{{ \rm  {sign}}}

\def\mod {{ \rm  {mod}}}

  \def\dim  {{ \rm  {dim}}}

\def\mod {{ \rm  {mod}}\,}

  \def\Tors {{ \rm  {Tors}}}

 \def\det {{ \rm  {det}}}

  \def\aug {{ \rm  {aug}}}

\begin{document}

\title{Torsions of 3-manifolds}
\author{Vladimir Turaev}
\address{Institut de Recherche Math\'ematique Avanc\'ee, Universit\'e 
Louis  
Pasteur - C.N.R.S.\\7 rue Ren\'e Descartes, F-67084 Strasbourg, 
France}
\asciiaddress{Institut de Recherche Mathematique Avancee, Universite 
Louis  
Pasteur - C.N.R.S.\\7 rue Rene Descartes, F-67084 Strasbourg, 
France}
\email{turaev@math.u-strasbg.fr}

\begin{abstract} We give a brief survey of  abelian 
torsions of 3-manifolds. \end{abstract}
\maketitle

\sh {Introduction}

This paper is a brief survey of my work on abelian torsions of 
3-manifolds. 
In 1976 I introduced an invariant, 
$\tau(M)$,  of a compact 
smooth (or PL-) manifold $M$ of any dimension 
(see \cite {Tu4} and references therein). This invariant is a sum of 
several Reidemeister torsions of $M$ numerated by   characters of the 
(finite 
abelian) group
$\Tors H$  where $H=H_1(M)$.   This invariant lies in a certain 
extension of 
the   group ring ${\bf {Z}} [H]$ and is defined   up 
to
multiplication  by $\pm 1$ and   elements of $H $. In the case $\dim 
M=3$ one can be more specific: if $b_1(M)\geq 2$ then 
 $\tau(M)\in {\bf {Z}} [H ]/\pm H$;    if $b_1(M)=0$ then 
 $\tau(M)\in {\bf {Q}} [H ]/\pm H$; if $b_1(M)=1$, then 
$\tau(M)$   can be expanded as a sum of 
an element of ${\bf {Z}} [H ]$ and a certain standard fraction. 
Classically, 
the Reidemeister torsions are used to 
distinguish homotopy equivalent but not simply homotopy equivalent 
spaces 
like lens spaces or their connected sums.  The study of    $\tau(M)$    
was motivated by its 
connections with the
Alexander-Fox
invariants of
$M$.  The 
present interest to this invariant is motivated by its connections to 
the 
Seiberg-Witten invariants.

To get rid of the ambiguity in the definition of  $\tau(M)$  
one needs to involve additional structures on $M$.  In  \cite {Tu1}  I 
introduced a
refined version 
$\tau(M,e, \omega)$ of $\tau(M)$ depending on the choice of  a 
so-called Euler structure $e$ on
$M$ and a homology orientation
 $\omega$ of 
$M$ (this  is an orientation in the real vector space $H_*(M;{\bf 
{R}})$). 
An Euler structure  on $M$ is a non-singular tangent vector field  on 
$M$ 
considered up to homotopy and an arbitrary modification in a small 
neighborhood of a point. The set $\Eul (M)$ of   Euler structures on 
$M$ has a natural involution $e\mapsto e^{-1}$ transforming the class 
of 
a non-singular   vector field to the class of the opposite  vector 
field. If $\chi(M)=0$, then the group $H=H_1(M)$ acts freely and 
transitively on 
 $\Eul (M)$ so that $\vert \Eul (M)\vert =\vert H\vert$. 

 The 
  invariant $\tau(M,e, \omega)$ has no indeterminacy  and   
$\tau(M)=\pm H \,
\tau(M,e, \omega)$ for all $e,
\omega$.  The torsions of various 
Euler 
structures on $M$ are computed from each other via $\tau(M,he, 
\omega)=h\, 
\tau(M,e, 
\omega)$ for any $h\in H$.

In this paper we shall   assume that $M$ is a closed connected 
oriented 
3-manifold. It has a   homology orientation $\omega_M$ 
defined by a   basis $ ([pt], b, b^*, [M])$ in\break $H_*(M; {\bf {R}})$
where $[pt]\in H_0(M; {\bf {R}})$ is the homology class of a point, $b$ 
is 
an 
arbitrary basis in $ H_1(M; {\bf {R}})$, $b^*$ is the  basis 
 in $H_2(M; {\bf {R}})$ dual to $b$      with respect to the 
(non-degenerate) intersection form $H_1(M; {\bf {R}})\times H_{2}(M; 
{\bf {R}})\to 
{\bf {R}}$, and finally $[M] \in H_3(M; {\bf {R}})$ is the fundamental 
class of 
$M$.
The homology orientation $\omega_M$ does not depend on the choice of 
$b$.
We shall write $\tau(M,e)$ for $\tau(M,e, \omega_M)$ where $e\in \Eul 
(M)$. If $-M$ is $M$ with opposite orientation, then 
$\tau(-M,e)=(-1)^{b_1(M)+1}\, \tau(M,e)$.

 The torsion $\tau(M,e)$ satisfies a fundamental duality formula 
 $$\overline {\tau(M,e)}=\tau (M,e^{-1})$$
 where   the overbar denotes the   conjugation in 
the 
group 
ring sending   group elements to their inverses.

The torsion $\tau (M)$ determines the first elementary ideal 
$E(\pi)\subset 
{\bf {Z}} [H]$ of the fundamental group $\pi=\pi_1(M)$. In particular 
if 
$b_1(M)\geq 1$, then  
$E(\pi)= \tau(M,e) I^2$ where $e$ is any Euler structure on $M$ and 
$I\subset  
{\bf {Z}} [H]$ 
is the  augmentation ideal. This 
implies that  $\tau(M,e)$  determines the   Alexander-Fox polynomial   
$\Delta(M)=\Delta(\pi)\in {\bf {Z}} [H/\Tors H]$. The torsion $\tau(M)$ 
can be viewed as a 
natural lift of   $\Delta(M)$ to ${\bf {Z}} [H]$. In contrast to 
$\Delta(M)$,
the torsion $\tau(M)$  in general is not determined by $\pi_1(M)$; 
indeed, it 
distinguishes lens spaces with the same $\pi_1$.

 The torsion $\tau(M)$    can be rewritten in terms of a numerical {\it 
torsion 
function} $T_M$ on the set $\Eul (M)$,
see 
\cite {Tu2}. This function takes values in ${\bf {Z}}$ if $b_1(M)\neq 
0$ and 
in ${\bf {Q}}$ if $b_1(M)= 0$. If $b_1(M)\neq 1$, then   $\tau$ and   
$T_M$ 
are related by 
the 
formula
$$\tau(M,e)=\sum_{h\in H} T_M(he)\, h^{-1} $$ 
for any $e\in \Eul (M)$. For $b_1(M)=1$,    there is a similar but a 
little more complicated formula. The   torsion 
function has a finite support and satisfies   the identity 
$T_M(e)= T_M(e^{-1})$ 
for all $e\in \Eul (M)$. 

For $b_1(M)=0$, we have  $\sum_{e\in \Eul (M)} T_M(e)=0$.  For 
$b_1(M)\geq 1$, 
the 
number  $\sum_{e\in \Eul (M)} T_M(e)$ is   essentially the 
Casson-Walker-Lescop  invariant $\lambda(M) \in {\bf {Q}}$:
$$ \sum_{e\in \Eul (M)} T_M(e)=\left \{ \begin {array} {ll} 
(-1)^{b_1(M)+1} \lambda (M) ,~ {\rm {if}} \,\,\, 
b_1(M)\geq 2, 
\\
\lambda 
(M)+ \vert \Tors H \vert/12
,~ 
{\rm
{if}} \,\,\, b_1(M)=1.  \end {array}  \right. $$

\sh {Connections to the Seiberg-Witten theory}

 The Seiberg-Witten 
invariant of 
$M$   is a numerical function, $ SW_M$,  on the set of  
$Spin^c$-structures on $M$, see for
instance     
 \cite {MT},  \cite {Lim},  \cite {MW}   and references therein.   For 
a 
$Spin^c$-structure 
$e$ on 
$M$, 
the 
integer  $SW_M(e )$ 
is the   algebraic    number
of solutions, called
monopoles,   to a certain system of differential equations associated 
with $e$. This number  
coincides with the  4-dimensional 
SW-invariant  of  the    $Spin^c$-structure $e \times 1$ on   $M\times 
S^1$. 

The invariants $ SW_M $ and
$\tau  (M)$  turn out to be 
   equivalent (at least up to sign).  The first step in this direction 
was made by Meng and Taubes
\cite {MT}  who observed that $SW_M$ determines the Alexander-Fox 
polynomial   $\Delta(M)$. The
equivalence between $SW_M$ 
and
$\tau(M)$ was established in \cite {Tu3}  in the case $b_1(M)\geq 1$. 
The Euler structures on $M$  are identified there with  
$Spin^c$-structures  on $M$ and  it is proven that $SW_M(e)=\pm T_M(e)$   
for 
all $e\in \Eul (M)$.
A similar
result for rational homology spheres was recently obtained by 
Nicolaescu \cite {Ni}.

For any $e\in Spin^c (M)$, the number $T_M (e)=\pm SW_M(e)$ can be 
viewed 
as 
the Euler
characteristic of the Seiberg-Witten-Floer monopole homology of $M$ 
associated 
with $e$ 
(see
\cite {MW}).  The same number appears also as the 
Euler
characteristic of the Floer-type homology of $M$ associated with $e$ by 
Ozsv\'ath 
and
Szab\'o, see \cite {OS}.

 \sh {Surgery formulas}

  The definition of $\tau $ is based on the   
methods of the theory of 
torsions, specifically,
cellular
chain complexes, coverings,  etc.  The definition of  $SW $ is 
analytical. These 
definitions     are  not
always suitable for explicit
computations.
We   outline   a  surgery formula  for  
$T_M$ suitable for   
computations. 

We first give a surgery description of Euler structures  (=  
$Spin^c$-structures) on   3-manifolds. To this end      
we introduce a notion of  a charge.
 A {\sl charge} 
on an oriented link
$L=L_1\cup ... \cup L_m\subset S^3$
is  an
$m$-tuple  $(k_1,...,k_m) \in {\bf   {Z}}^m$ such that for all
$i=1,...,m$,  
$$k_i\equiv 1+ \sum_{j\neq i} lk (L_i, L_j)\, (\mod 2)  $$
where $lk$ is the linking number in $S^3$. 
A charge $k$ on $L$    determines  an Euler structure,  
$e_k^M$, on any 3-manifold
$M$  obtained by surgery 
on         
$L$, see \cite {Tu5}. 
 
The surgery  formula computes
$T_M(e_k^M)$   in terms of the framing and  linking 
numbers of the components of
$L$, and the
Alexander-Conway polynomials of $L$ and its sublinks.  
Thus, the algebraic number of monopoles can be computed (at least up to 
sign) in terms of
classical  link invariants. For simplicity, we state here the surgery 
formula only 
in the case of algebraically split links and the first Betti number 
$\geq 2$.

Let $L=L_1\cup ... \cup L_m\subset S^3$ be an oriented algebraically 
split 
link 
(i.e.,    $lk(L_i,L_j)=0$ for all 
$i\neq j$).
  Recall    the
Alexander-Conway polynomial  $\nabla_L  
 \in {\bf   {Z}} [t_1^{\pm 1},..., t_m^{\pm 1}]$, see \cite {Ha}.
Since  $L$ is  algebraically split,  
$\nabla_L $ is   divisible by $ \prod_{i=1}^m (t_i^2-1)$ in ${\bf   
{Z}} 
[t_1^{\pm 1},..., t_m^{\pm
1}]$. 
 We have a  finite  expansion  
 $$\nabla_L / \prod_{i=1}^m (t_i^2-1)=\sum_{l=(l_1,...,l_m)\in {\bf   
{Z}}^m} 
z_l (L) \,\, 
t_1^{l_1}... t_m^{l_m}$$
where $z_l (L)\in {\bf   {Z}}$. 
 
Let    $M$ be    
obtained  by surgery on $L$ with framing    
$f=(f_1,...,f_m)\in {\bf   {Z}}^m$.   
Denote by
$J_0$  the  set of  all $j\in \{1,...,m\}$   such that  
$f_j=0$. For a   set $J\subset \{1,...,m\}$,  denote the link 
 $\cup_{j\in  J} L_j$ by $L^J$.  Put     $\vert J 
\vert=\card (J)$ and suppose that $\vert J_0\vert \geq 2$.
Then for any charge  $k=(k_1,...,k_m)$ on $L$, 
$$  T_M(e_k^M )\kern 4.3in \relax\hbox{}   \eqno (1) $$
$$= (-1)^{m+1}  \sum_{J_0\subset J \subset \{1,...,m\}}\,\, (-1)^{\vert 
J 
\backslash J_0  
\vert}\, \prod_{j\in   J\backslash J_0} 
\sign (f_j)  \,\,  
  \sum_{l\in {\bf   {Z}}^{J}, l= -k\, (\mod 2f)} z_{l} (L^{J}). $$ 
 Here the sum goes over  
 all sets   $J\subset \{1,...,m\}$ containing $J_0$. The sign   $\sign 
(f_j)=\pm   1$   of $f_j$  
is
well defined  since $f_j\neq 0$ for $j\in J\backslash J_0$. The formula
$l\in {\bf   {Z}}^{ J}, l= -k \,(\mod 2f)$   means that  
$l$ runs over all tuples of integers numerated by elements of $J$   
such that $l_j= -k_j \,
(\mod 2f_j)$ for all
$j\in J$.   By
$\vert J  \vert\geq \vert J_0  \vert \geq 2$, the algebraically 
split link $L^{J}$  has 
$\geq
2$ components  so that 
$z_{l} (L^{J})$ is a well defined  integer.  Only a finite number of 
these integers  are non-zero
and therefore
the sum in (1) is finite. This sum  obviously 
depends only on $k (\mod
2f)$; the Euler structure 
 $e_k^M$ also depends  only on $k (\mod 2f)$.

For   a link  $L =L_1\cup ... \cup L_m$ which is not algebraically 
split,  the polynomial
$\nabla_L $ can be   divided by $
\prod_{i=1}^m (t_i^2-1)$ in a certain quotient  of   ${\bf   {Z}} 
[t_1^{\pm 
1},..., t_m^{\pm 1}]$. 
This leads to a    surgery
formula for  an arbitrary $L$, see \cite {Tu5}.
Formula (1) and its generalizations to non-algebraically split links 
yield 
similar formulas for the   Alexander polynomial  
$\Delta(M)$ and  the Casson-Walker-Lescop invariant of $M$ (in the case 
$b_1(M)\neq 0$).

\sh{Moments of the torsion function}

Every $e\in \Eul (M)$ has a characteristic class $c(e)\in 
H=H_1(M)$ defined as the unique element of $H$ such that $e=c(e) 
e^{-1}$.  
This class is the first obstruction to the existence 
of a homotopy between a vector field   representing $e$ and the 
opposite vector field.
For
any $x_1,..., x_m\in H^1(M;\bf   R)$, we define the corresponding 
{\it $m$-th moment} 
of
$T_M$ by $$\langle T_M \, \vert \, x_1,..., x_m \rangle= \sum_{e\in 
\Eul (M) } \, 
T_M ( e)
\prod_{i=1}^m \langle c(e), x_i \rangle  $$ where on the 
right-hand 
side
$\langle \,, \, \rangle$ is the evaluation pairing $H  \times H^1 
(M;{\bf  {R}})\to \bf   R$.  
It turns out  that if $ m\leq b_1(M)-4$, then $\langle T_M \, \vert \, 
x_1,..., x_m
\rangle=0$.  In particular, if $b_1(M)\geq 4$ then $\sum_{e } \, T_M ( 
e)=0$.

  Interesting phenomena occur for $ m= b_1(M)-3$.  Set $b=b_1(M)$.  If 
$b$ is 
even
then $\langle T_M \, \vert \,x_1,..., x_{b-3} \rangle=0$ for any 
$x_1,...,
x_{b-3}\in H^1(M;\bf   R)$.  For odd $ b$, the number $\langle T_M \, 
\vert \,
x_1,..., x_{b-3} \rangle$ is determined by the cohomology ring of 
$M$ with
coefficients in ${\bf   {Z}}$.  We state here a special case of this 
computation.
Recall that an element of $H^1(M)$ is {\sl primitive} if it is 
divisible only by
$\pm 1$.  If $b\geq 3$ is odd then for any primitive $x\in H^1(M)$, 
$$\langle T_M
\, \vert \, \underbrace {x ,..., x}_{b-3} \rangle =  2^{b-3} (b-3)!  
\,\vert 
\Tors H \vert \, \det 
\,g_x 
\eqno
(2)$$ where $g_x$ is the skew-symmetric bilinear form on the lattice
$H^1(M)/{\bf   {Z}} x$ defined by $g_x (y,z)= (x\cup y \cup z) ([M])$ 
for 
$y,z\in
H^1(M)/{\bf   {Z}} x$.  Note that $\det \,g_x= (\Pf (g_x))^2\geq 0$.  
This 
computation
implies for instance that if $x$ is dual to the homology class of a 
(singular)
closed oriented surface in $M$ of genus $ \leq (b-3)/2$, then $\sum_{e 
} \, T_M (
e) \langle c(e), x \rangle^{b-3} =0$.  For $b=3$, formula (2) gives 
$$\sum_{e }
\, T_M ( e)= \vert \Tors H \vert \, ((x\cup y \cup z) ([M]))^2 \eqno 
(3)$$ 
where 
$x,y,z$
is any basis of $H^1(M)$. Formula (3) and  the equality 
$\sum_{e }
\, T_M ( e)=\lambda(M)$ yield Lescop's computation of $\lambda(M)$ for 
$b_1(M)=3$.

\sh{Basic Euler structures and the Thurston norm}

  An Euler structure 
$e\in \Eul (M)$ is said to be {\sl basic} if $T_M(e)\neq 
0$. The set of basic Euler structures is closely related to the 
Thurston norm, see \cite {Th}.    Recall that the Thurston norm of 
$s\in 
H^1(M)$ 
is 
defined 
by $$\vert\vert s
\vert\vert_T=\min_S \,\{\chi_-(S)\}$$ where $S$ runs over  
closed oriented
embedded (not necessarily connected) 
surfaces in $ M$   dual to $s$ and $\chi_- 
(S)=\sum_i \max (-\chi 
(S_i), 
0) 
$ where 
the
sum runs over all components $S_i$ of $S$.  Then
  for any
$s\in H^1(M)$ and any basic Euler structure $e$ on $M$, 
$$\vert\vert s \vert\vert_T \geq  \vert \langle c(e), s 
\rangle \vert. \eqno (4)$$ This inequality is a
cousin of the classical Seifert inequality    which says 
that the genus of a knot in $  S^3$ is greater 
than or 
equal
to the half of the span of its Alexander polynomial. The inequality (4) 
is a 
3-dimensional version of the (much deeper) adjunction inequality in 
dimension 
4. 
A weaker version of (4)  involving    $\Delta (M)$ rather 
than 
$\tau(M)$  appeared in \cite {McM}.  
 For  more general homological estimates of the Thurston norm,
see \cite {Har}, \cite {Tu6}. 
For analogous     
estimates 
  in the Seiberg-Witten theory in dimension 3,  
 see 
\cite {Au},  \cite {Kr}, \cite {KM1},
\cite {KM2}. Similar  estimates  appear also in the   Ozsv\'ath-Szab\'o  
theory 
\cite {OS}.

\sh{Examples}

Let 
$M$ be 
the
total space of an oriented circle bundle over a closed connected 
orientable
surface $\Sigma$ of genus $g\geq 0$. Let $e_{\pm} \in  \Eul (M)$ be 
represented 
by the non-singular 
vector 
field  
on $M$ tangent   to the fibers of the bundle $M\to \Sigma$ in the 
positive 
(resp. 
negative) direction. Observe that $e_{-}=(e_{+})^{-1}$ and
$c(e_-)=e_-/e_+=t^{ 2g-2}$ where $t \in H =H_1(M) $ is 
the 
homology
class of the   fiber $S^1$.   We claim that $$ \pm \tau (M,e_-)=
  \pm  (t-1)^{2g-2}. \eqno (5) $$  Here we do not (homologically) 
orient $M$ 
and 
consider 
the torsion  up to sign. Applying (5) to the   opposite orientation 
of  the 
fibers, 
we 
obtain that $\pm \tau (M,e_+)$ $=
  \pm (t^{-1}-1)^{2g-2} $. The same formula follows from (5) and the 
duality 
for $\tau$. 

 The Thurston norm for $M$ can be easily computed  since most (if not 
all) 
generators of $H_2(M)$ are represented by tori. For any $s\in 
H_1(M;{\bf {R}})$, we have 
$\vert\vert s
\vert\vert_T=\vert \langle t, s\rangle \vert$. In particular   if the 
bundle 
$M\to \Sigma$ is non-trivial then the Thurston norm  is identically 0.
As an exercise, the reader may compute the torsion function for $M$ (at 
least 
up to sign) and check (4). Similar computations are available in the   
Seiberg-Witten theory, see \cite {Ba}.
 
In particular, if $M$ is the 3-torus $S^1 \times S^1 \times S^1$ and $e 
\in  
\Eul (M)$ is 
represented 
by the non-singular 
vector 
field  
on $M$ tangent   to the fibers of the obvious projection $M\to  S^1 
\times 
S^1$, then $\pm \tau(M,e)=\pm 1$. This and formula (3) imply that 
$\tau(M,e)= 
 1$ for any orientation of $M$.

\sh{Realization}

The realization problem for $\tau$ consists in finding necessary and 
sufficient conditions for a pair (a finitely 
generated
abelian group, an element   of its group ring) to be realizable as the   
first homology group and the torsion $\tau$  of a   closed connected 
oriented 3-manifold.  In spirit of Levine's \cite {Lev} realization 
theorem  for 
the Alexander polynomial  of 
links  in $S^3$, 
we have the following partial
 result.

{\it Let $H$ be a 
finitely 
generated
abelian group   and $\lambda\in {\bf {Z}}[H]$ be symmetric 
with  
$\aug(\lambda)=1$.
  If a pair $(H,\tau )$ is realizible   then so is  
$(H, \lambda\,
\tau) $.}

Here $\lambda$ is said to be {\it symmetric} if  $\overline 
\lambda=\lambda$.  For example, since 
the torsion of a 3-torus is $1$, we obtain that   any  
symmetric $\lambda\in
{\bf  {Z}}[{\bf  {Z}}^3]$ with $\aug(\lambda)=1$   is 
realizable   as the torsion $\tau$  of a closed oriented 3-manifold 
with 
$H_1 ={\bf  {Z}}^3$. 
It follows from the Bailey theorem (see \cite {Hi}) and the surgery 
formula 
outlined above that 
 any symmetric  
$\lambda\in
{\bf  {Z}}[{\bf  {Z}}^2]$   is 
realizable   as the torsion $\tau$  of a closed oriented 3-manifold 
with 
$H_1 ={\bf  {Z}}^2$.

\Addresses\recd
 
 \end{document}